\documentclass[12pt]{amsart}

\usepackage{amsmath,amssymb,amsthm}
\usepackage{mathrsfs, wrapfig, graphicx, caption, cite, pifont}
\usepackage{paralist}
\usepackage[all]{xy}

\numberwithin{equation}{section}

\newtheorem{theorem}{Theorem}[section]

\newtheorem{lemma}[theorem]{Lemma}

\theoremstyle{definition}
\newtheorem{definition}[theorem]{Definition}

\theoremstyle{remark}
\newtheorem{rem}[theorem]{Remark}

\newtheorem{exam}[theorem]{Example}

\newcommand{\be}{\begin{equation}}
\newcommand{\ee}{\end{equation}}

\def\beq{\begin{equation}
}
\def\eeq{\end{equation}}

\def\supp{\hbox{supp}}

\def\dl{\delta}
\def\idl{\frac{1}{\delta}}

\def\eql{equilibrium}
\def\qed{$\blacksquare$}

\def\pt{\partial}

\def\RR{\mathbb R}

\def\ben{\begin{enumerate} }

\def\een{\end{enumerate} }

\def\RR{ {\mathbb R} }
\def\Rn{ {\mathbb R}^n }
\def\Rnp{ {\mathbb R}^n_{>0} }
\def\Rnn{ {\mathbb R}^n_{\geq 0} }

\def\bRplus{ {\mathbb R}^n_{ >0} }
\def\bRplusc{ {\mathbb R}^n_{ \geq 0} }

\def\lout{\Lambda_{o}}

\begin{document}

\title[Homotopy and counting
equilibria in reaction networks]
{Homotopy methods for counting
reaction network equilibria}

\author[Craciun]{Gheorghe Craciun}
\address{Department of Mathematics and Department of Biomolecular
Chemistry, University of Wisconsin-Madison.
(craciun@math.wisc.edu)}

\author[Helton]{J. William Helton}
\address{Mathematics Department, University of California at San
Diego, La Jolla CA 92093-0112
(helton@ucsd.edu)}

\author[Williams]{Ruth J. Williams}
\address{Mathematics Department, University of California at San
Diego, La Jolla CA 92093-0112
(williams@euclid.ucsd.edu)}

\maketitle

\centerline{\large \today}

\begin{abstract}
Dynamical system models of complex biochemical reaction networks are usually
high-dimensional, nonlinear, and contain many unknown parameters.
 In some cases the
reaction network structure dictates that positive equilibria must be unique for all
values of the parameters in the model.
In other cases multiple equilibria exist if and
only if special relationships between these parameters are satisfied.
We describe
methods based on homotopy invariance of degree which allow us
to determine the number of
 equilibria for complex biochemical reaction networks
and how this number depends on
parameters in the model.
\end{abstract}

\vfill\break
\tableofcontents

\vfill\break

\section{Introduction}
\label{sec:intro}

Dynamical system models of complex biochemical reaction networks are usually
high-dimensional, nonlinear, and contain many unknown parameters.
As was shown recently in \cite{CTF06}, based on the assumption
of mass-action kinetics, graph-theoretical properties of some
biochemical reaction networks can guarantee the uniqueness of
positive equilibrium points for any values of the reaction rate
parameters in the model.
On the other hand, relatively simple
reaction networks do admit multiple positive equilibria for some
values of the parameters
as shown in  \cite{CF05, CF06, CTF06}.

The aforementioned results do not address the dependence of the number of
equilibria on the parameter values unless there
is  a unique equilibrium for every set of parameters.
Also they  do not address the general problem of \it existence \rm
of positive equilibria.
Here we describe
methods using  degree theory
to analyze  general biochemical dynamics
(not only mass-action kinetics).
These methods allow us to determine how the number of
positive equilibria for a complex biochemical reaction
network depends on the parameters of the model.
They will often also imply the existence of positive
equilibria.
Also we obtain uniqueness  of positive equilibria
in various situations under significantly
 weaker assumptions than in  \cite{CF05, CF06, CTF06}.

\subsection{Overview}

We are interested in equilibria
for high-dimensional, nonlinear dynamical systems
that originate from chemical dynamics.
These dynamical systems are systems of
ordinary differential equations of the form
\beq
\label{eq:diffeq}
\frac{dc}{dt}= f(c) \quad
\eeq
where
$f$ is a smooth function defined on a subset
of the orthant $\mathbb{R}^n_{\geq 0}$
of vectors $c$ in $\mathbb{R}^n$ having nonnegative components.
Such dynamical systems usually have a large number of state variables,
i.e., $n$ is large.
In addition, the parameters defining $f$ are often not well known.
The focus of this paper is on equilibria for dynamical systems of the form
(\ref{eq:diffeq}), that
is on $c^\ast $ for which $f(c^\ast) =0$.

We consider the  dynamical system  (\ref{eq:diffeq}) on a subset
$\overline\Omega$ of $\Rnn$ which is the closure of
a domain $\Omega$
in $\Rnp$.
We give conditions for the number of equilibria of (\ref{eq:diffeq}) to
remain constant 
as we ``continuously
deform" (homotopy) the function $f$ through a family
of functions.
A key assumption is
that
the following condition  holds for all members of the family:
\begin{quote}
{\bf(DetSign)}
\ \ \ {\it The determinant of the Jacobian matrix
$\frac{\pt f}{\pt c}(\cdot)$ of $f$
is either strictly positive or strictly negative on  $\Omega$.
%
}
\end{quote}
(Recall that the Jacobian $\frac{\pt f}{\pt c}(c)$ at $c$
is the matrix
$ \{ \frac{\pt f_j}{\pt c_i}(c), \  i,j=1,\dots, n\} $.)

\it
What we observe is that when the condition
(DetSign) holds for all $f$ in the family and $\Omega $ is bounded,
then the number of equilibria for the dynamical system
(\ref{eq:diffeq})
is a constant for all $f$ in the family,
provided
there are no equilibria on the boundary of
$\Omega$ for any  $f$ in the family (see Theorem \ref{thm:intro}). 
We further indicate
how this result extends to unbounded domains such as $\Omega =\Rnp$ under
suitable  conditions (see e.g., Theorem \ref{thm:dissip}),
including those associated with
a mass-conserving reaction network operating in a chemical
reactor with inflows and outflows. \rm

This paper extends
previous findings in several ways.

The (DetSign) condition  was introduced  by
Craciun and Feinberg \cite{CF05, CF06} in the
context of chemistry with $\Omega=\Rnp$ and they
observed  that
many chemical reaction networks have the property
(DetSign).
They gave many examples and many tests
for this condition to hold in the case where $f$ is a system of polynomials
and $\Omega$ is  $\Rnp$.
Then they
\cite{CF05, CF06} proved that if
the components of $f$ are polynomials corresponding to
mass-action  kinetics (operating in
a continuous flow stirred tank reactor), and
if (DetSign) holds on $\Rnp$ for \it all positive \rm ``rate constants",
then for each particular choice of  rate constant,
when an  equilibrium exists, it is unique.
Here we obtain  stronger conclusions
with  weaker assumptions. In particular the following are
features of our approach.
\ben
\item
Rather than
all positive ``rate constants"
we can select a (vector valued) rate constant $k_0$ of interest
at which (DetSign) holds.
Then one merely needs
a continuous curve  $k(\lambda)$
of ``rate constants" joining $k_0$
to another $k_1$ at which (DetSign) holds
and at which the dynamical system has
a unique positive equilibrium.
\item
For a mass-conserving reaction network operating
in a chemical reactor with
inflows and outflows, under the (DetSign) assumption in (1),
we prove existence
and uniqueness of a positive equilibrium, see Theorem \ref{thm:mm}.
\item
In (1) and (2),  the function
$f$ need not be polynomial, but is required
only to be continuously differentiable.
Of practical importance are rational $f$
as one finds in Michaelis-Menten or Hill type
chemical models,
see \S \ref{sec:dynam}, \S \ref{sec:appl}.
\item
We give methods, see \S \ref{sec:appl}, combining
the items above to describe large regions of rate
constants where a chemical  reaction network has
a unique positive equilibrium.
\een

We also point out in this paper that the biochemical reaction network
  models introduced
and analyzed by Arcak and Sontag  \cite{arcson06, ASprept}
satisfy (DetSign) and we can also rule out
boundary equilibria (where they give enough data).
Consequently,  under extremely weak
hypotheses, we obtain that each of these models has a unique
positive equilibrium,
see \S \ref{sec:ASsummary}.
The findings of  Arcak and Sontag are impressive
in that under strong hypotheses  they prove
global asymptotic stability of equilibria, a topic that this paper
does not address.


\subsection{More detail}

Now we give some formal definitions.
Let $\Omega$ be a \bf domain  \rm in
$\Rn$, i.e., an open, connected set in $\Rn$.
We denote the closure of $\Omega$ by $\overline \Omega$
and the boundary of $\Omega $ by $\partial \Omega$.
A function $f:\overline \Omega \to \Rn$ is \bf smooth \rm
if it is once continuously differentiable on $\overline \Omega$.
If $\Omega $ is bounded, for such a smooth function $f$,
the following norm is finite:
$$
\| f \|_\Omega:= \sup_{c \in \Omega} \|f(c) \|
$$
Here $\| \cdot \|$ denotes the Euclidean norm on $\RR^n$.
When $\Omega $ is bounded,  a family
$f_\lambda : \overline \Omega \to {\mathbb R}^n$ for
$\lambda \in [0,1]$,
is  a    {\bf continuously varying family of functions}
provided
each $f_\lambda$ is smooth
and the mapping $\lambda\to f_\lambda$  is continuous
on $[0,1]$ with the norm  $\| \  \cdot \|_\Omega$ on the functions
$f_\lambda$.
A \bf zero \rm of $f:\overline \Omega \to \Rn$ is a value
$c\in \overline \Omega$ such that $f(c) =0$, where $0$ is the
zero vector in $\Rn$.
A zero of $f$ is an equilibrium point for the dynamical system
(\ref{eq:diffeq}).

The following is an immediate consequence of Theorem
\ref{thm:degreeDet} which will be proved in \S \ref{sec:degree}.
This theorem and
examples  given in \S \ref{sec:exSont} are designed to
illustrate our approach;
then more targeted theorems are given in \S \ref{sec:mass} and
\S \ref{sec:dynam}, followed by more examples in \S \ref{sec:appl}.

\begin{theorem}
\label{thm:intro}
Suppose $\Omega \subset {\mathbb R}^n$
is a bounded domain and
$f_\lambda : \overline\Omega \to {\mathbb R}^n$ for
$\lambda \in [0,1]$,
is a continuously varying family of smooth functions
such that $f_\lambda$ does not have any zeros on
the boundary of $\Omega$ for all $\lambda\in [0,1]$.
If $\det \left( \frac{\partial f_\lambda}{\partial c} (c) \right) \neq 0$
for all $c \in \Omega$
whenever $\lambda=0$ and whenever  $\lambda=1$,
then the number of zeros of $f_0$ in $\Omega$ equals
the number of zeros of $f_1$ in $\Omega$.
\end{theorem}

The domain of interest for chemical dynamics (cf. (\ref{eq:diffeq}))
is typically the orthant $\Rnp$,
but this is not bounded, so it violates the
hypothesis ``$\Omega$ is bounded".
Thus in applying Theorem
\ref{thm:intro} we must approximate $\Rnp$
by a large bounded domain $\Omega$
and check for the absence of boundary equilibria.
One can think of the boundary  $\partial \Omega$
in two pieces:
that which intersects the boundary of $\Rnp$, called the {\bf sides
of $\Omega$},
and the {\bf outer boundary},
$\partial \Omega \cap  \Rnp $.
We show that if we assume conservation of mass (e.g., by atomic balance)
in our model and augment with suitable ``outflows",
then natural bounded domains $\Omega$ can be chosen
which have no equilibria on the outer boundary.
An example of such a natural bounded domain in $\mathbb{R}_{>0}^2$
is shown in Figure \ref{bdom}.
Also under assumptions of positive invariance and augmentation with
``inflows",
there are no equilibria on the sides.  In these cases we conclude
that there is exactly one nonnegative equilibrium $c^\ast$ for (\ref{eq:diffeq})
and that
it is actually a positive equilibrium,
i.e., it lies  in $\mathbb{R}_{> 0}^n$.
This result is described in detail
in \S \ref{sec:mass} and \S \ref{sec:dynam}.

{\unitlength=0.3mm 
\begin{figure}[htbp]
\begin{picture}(400,200)(10,20)
\put(135,55){\vector(1,0){230}}
\put(135,55){\vector(0,1){135}}
\put(135,155){\line(2,-1){200}}
\put(255,128){\makebox(0,0){{$\longleftarrow$ Outer boundary}}}
\put(230,45){\makebox(0,0){$\big\uparrow$}}
\put(230,30){\makebox(0,0){{A side of $\Omega$}}}
\put(350,40){\makebox(0,0){{$c_1$}}}
\put(85,105){\makebox(0,0){{A side of $\Omega \ \longrightarrow$}}}
\put(110,180){\makebox(0,0){{$c_2$}}}
\put(190,90){\makebox(0,0){{$\Omega$}}}
\end{picture}
\caption{An example of a natural bounded domain in $\mathbb{R}^2_{>0}$ with the outer boundary and sides indicated.}
\label{bdom}
\end{figure}
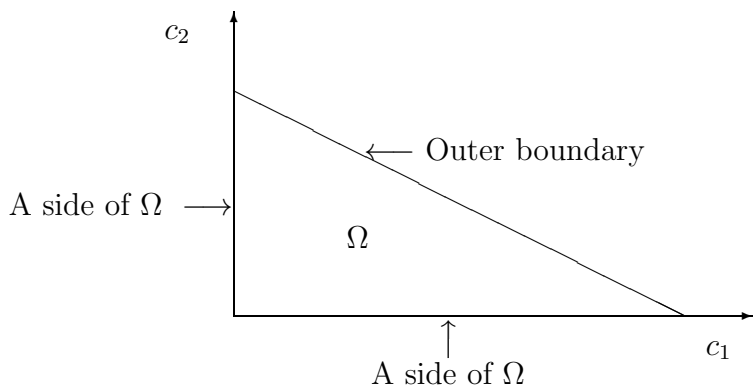}

\subsection{Organization of the paper}

In this paper we give many examples of widely varying types
to illustrate our contention that our method applies broadly and is easy
to use. Section \ref{sec:exSont}  gives several examples
illustrating Theorem \ref{thm:intro}.
Section   \ref{sec:degree}  summarizes degree theory
since our proof is based on this
and relies on the observation that when (DetSign) is true, and there
are no boundary equilibria, then
``the degree of $f$ with respect to 0"
equals $\pm$ the number of equilibria in a bounded domain.
Then in \S \ref{sec:degree} we prove Theorem \ref{thm:intro} and more.
Section \ref{sec:mass}
describes the mathematical benefits of mass dissipation (including
mass conservation)
and ``inflows" and ``outflows".
Section  \ref{sec:dynam} describes a chemical
reaction network framework which contains
in addition to
mass-action kinetics,  Michaelis-Menten and Hill  dynamics.
We conclude in \S  \ref{sec:appl} with more examples and  new methods
presented in the context of these examples.
For many of our examples, the determinant
of the Jacobian of $f$ was computed symbolically using Mathematica.
As a complement to this paper, we have established a webpage at 
http://www.math.ucsd.edu/$\sim$helton/chemjac.html
containing Mathematica notebooks
for  many examples in this paper and a demonstration notebook that
readers may edit to run their own examples.


\section{Examples}
\label{sec:exSont}

Our goal in this section is to present some examples
showing how to use Theorem \ref{thm:intro}.
In the process we mention that all chemical reaction
examples of Arcak and Sontag \cite{arcson06, ASprept}
satisfy (DetSign)  and
fit well into our approach here.
Later in \S \ref{sec:appl} we give broader categories of
examples.

\subsection{Two examples on
treating boundary equilibria}

We start with two examples, the study of which goes back
to a class of examples studied by Thron \cite{tysoth78, thr91}.
Here, $c$ satisfies (\ref{eq:diffeq})
and
the Jacobian  for all $c$ has the form:
\begin{equation}
     \label{eq:Amatrix}
       \frac{\partial f}{\partial c} \; = \;
       \left[\begin{array}{ccccc} -a_1  & 0 & \cdots & 0 & -b_n
      \\
b_1 & -a_2 &\ddots & & 0 \\
       0 & b_2 & -a_3 & \ddots & \vdots \\ \vdots & \ddots & \ddots &
       \ddots & 0 \\ 0 & \cdots & 0 & b_{n-1} & -a_n
       \end{array} \right]
       \end{equation}
where $a_i\geq 0, $
       $b_i \geq 0,$
       $i=1,\ldots,n,$ may depend on $c$.
This cyclic feedback structure is common in gene regulation networks,
cellular signaling pathways, and metabolic pathways \cite{ASprept}.
Thron  showed that all eigenvalues of $\frac{\partial f}{\partial c}$
       have nonnegative real part (local stability)  if
       $ \frac{b_1\cdots b_n}{a_1\cdots a_n} \, < \, (\sec(\pi/n))^n.$
 Arcak and Sontag
        showed that an equilibrium of  such a  dynamical system
 is unique and globally stable under strong global restrictions.

       Here we observe that our key assumption,
       the determinant of the Jacobian never changes sign, is met,
which is a major step toward checking when a unique positive
equilibrium exists for this class of problem.

\begin{lemma} When the Jacobian has the form (\ref{eq:Amatrix}), we have
\be
\label{eq:detAS}
    \det \; \left( \frac{\partial f}{\partial c}\right)
    = \ (-1)^n \; [ \Pi_{j=1}^n a_j + \Pi_{j=1}^n b_j ]
\ee
which if not zero has
  sign  independent of $a_i, b_i \geq 0$.
\end{lemma}
{\bf Proof:} \
Direct computation.
\qed

We
 shall now consider several examples from papers of Arcak and Sontag
 primarily to illustrate
that  checking for absence of boundary equilibria is straightforward;
subsequently we obtain existence and uniqueness of an equilibrium.
In this section we assume that
all parameters in the reactions are strictly positive.
In the sequel, we shall
use $\dot c$ as an abbreviation for the time derivative of $c$.

\begin{exam}
\label{ex:AS1}
   \rm
We start with an example from \S6 of \cite{arcson06}
which they took from  Thron \cite{thr91}. For this,
\begin{eqnarray}
\dot{c}_1&=&  \frac{p_1c_0}{p_2 +c_3} -p_3 c_1 \label{mapk1}\\
\dot{c}_2&=&   p_3c_1  - p_4 c_2 \label{mapk2}\\
\dot{c}_3&=& p_4 c_2 - \frac{p_5 c_3}{p_ 6 + c_3}.
\label{mapk3}
\end{eqnarray}
Now we apply Theorem \ref{thm:intro}
    to obtain the conclusion:

\begin{quote}
{\it
For each set of parameters $p_j >0, j=1,2,\ldots,6 , c_0>0 $,
there is a
unique equilibrium
point $c^\ast$ in
$\mathbb{R}^3_{>0}$ for the chemical reaction
network with dynamics given by (\ref{mapk1})--(\ref{mapk3}),
and there is no equilibrium point
on the boundary of $\mathbb{R}^3_{\geq 0}$.
}
\end{quote}

However, we have made no comment on stability (even local),
while \cite{arcson06} gives  certain conditions ensuring  global
stability. (In fact, one can algebraically solve for the two
equilibria of (\ref{mapk1})--(\ref{mapk3})
as functions of the parameters.
Inspection reveals that exactly one of these is in $\mathbb{R}^3_{>0}$
and neither is on the boundary of $\mathbb{R}^3_{\geq 0}$.
The point
of us treating this example is to show in a simple context how our
method works.)

\noindent
{\bf Proof:} \ \
%
%
We shall apply Theorem \ref{thm:intro} to prove
that there is a unique equilibrium, inside any sufficiently
large box; hence there is a  unique equilibrium in the orthant.
The right hand side $f_{p, c_0}(c)$ of the
differential equations
(\ref{mapk1})--(\ref{mapk3}),
while a function of $c$, also depends on positive parameters
$(p,c_0)$.
One can check  that the Jacobian for any of these parameters
has the form
in (\ref{eq:Amatrix}) for all $c\in \mathbb{R}^3_{>0}$, and
since $n=3$ and all parameters are strictly positive,
the Jacobian  determinant is strictly negative.
Note that for any two values of the positive parameters,
$(p^*,c^*_0)$ and $(p^\dag,c^\dag_0)$,
$f_{\lambda(p^*,c_0^*)+(1-\lambda)(p^\dag,c^\dag_0)}, \lambda\in [0,1],$ defines a continuously
varying family of smooth functions on any bounded subset of $\Rnn$.
We check below that the no equilibria 
(i.e., no zeros of $f_{(p,c_0)}$) on the boundary hypothesis
holds  on any sufficiently large box, for all positive parameters $(p,c_0)$,
and thereby conclude using Theorem \ref{thm:intro}
that
the number of equilibria of (\ref{mapk1})--(\ref{mapk3}) in $\Rnp$ does not depend on $(p,c_0)$
provided the parameters are all strictly positive.
Computing the equilibria at one
simple ``initial" value of $(p^*,c^*_0)$ then finishes the proof.

\noindent
{\it No equilibria on the boundary of the orthant:}
    Suppose an equilibrium has $c_2=0$.
Then equation (\ref{mapk2}) implies $c_1=0$ which contradicts
(\ref{mapk1}).
Likewise if we start by assuming $c_1=0$ we get $c_1>0$
and a contradiction.
On the other hand if $c_3=0$,
    then (\ref{mapk3}) implies $c_2=0$, which reverts to the
case
considered first.  Thus there are no equilibria on the
boundary of $\mathbb{R}^3_{\geq 0}$.

\noindent
{\it No equilibria on the outer boundary of some big box:}
Suppose $0< \delta< \frac12$ and $\idl >p_j > \dl$ for all $j$ and $c_0 <1$.
Pick $\Omega$ to be a box
$\Omega:= \{ c \in \mathbb{R}^3_{>0}: c_j
< (\frac{1}{\delta})^4$ for $j=1,2,3\}$.
An \eql \ on the outer boundary of the box satisfies
\begin{enumerate}
     \item
$c_1=  (\idl)^4$ which by (\ref{mapk1}) implies
\ $ (\idl)^2 > \frac{p_1c_0}{p_2 +c_3}= p_3 c_1> (\idl)^3.$
A contradiction.\\ OR
\item
$c_2=  (\idl)^4$ which by (\ref{mapk3}) implies
\ $ \idl  >  \frac{p_5 c_3}{p_ 6 + c_3} =p_4 c_2 > (\idl)^3 .$
A contradiction.\\ OR
\item
$c_3=  (\idl)^4$ which by adding
(\ref{mapk1}), (\ref{mapk2}), (\ref{mapk3}) implies
that
$$\delta^3= \frac{\frac{1}{\delta}}{\frac{1}{\delta^4}}
\geq
\frac{p_1 c_0}{p_2 +c_3}
= \frac{p_5 c_3}{p_6 + c_3}
\geq
\frac{\delta \cdot \frac{1}{\delta^4}}{\frac{1}{\delta} +
\frac{1}{\delta^4}}
=\frac{\delta}{\delta^3 +1}.$$
A contradiction.
\end{enumerate}


\noindent
{\it Initializing:}
Up to this point, Theorem \ref{thm:intro}  tells us that
each choice of parameters yields the same number of
equilibria!  It is easy to compute for oneself that there is a simple
choice of parameters which yields a unique positive equilibrium,
for example $p_j=1$ for all $j$ yields the unique equilibrium,
$c_1=c_2 =\frac{c_0}{1+c_0}, c_3=c_0$.
Thus there is one and only one equilibrium
in $\mathbb{R}^3_{>0}$ for each value of the positive parameters $(p,c_0)$.
\qed
\end{exam}

\begin{exam} \ \
\label{mapk}
In Example 1 of \S 4 in \cite{ASprept},
\rm
the authors describe  a simplified model
of mitogen activated protein kinase (MAPK) cascades with inhibitory
feedback, proposed in \cite{KholodenkoNegfeedback, stas2}. For this,
\begin{eqnarray}
\dot{c}_1&=&-\frac{b_1c_1}{c_1+a_1}+\frac{d_1(1-c_1)}{e_1+(1-c_1)}
\frac{\mu}{1+kc_3} \label{amapk1}\\
\dot{c}_2&=&-\frac{b_2c_2}{c_2+a_2}+\frac{d_2(1-c_2)}
{e_2+(1-c_2)}c_1 \label{amapk2}\\
\dot{c}_3&=&-\frac{b_3c_3}{c_3+a_3}+\frac{d_3(1-c_3)}{e_3+(1-c_3)}c_2.
\label{amapk3}
\end{eqnarray}
The variables $c_j \in [0,1], \ j=1,2,3$ denote the
concentrations of the active forms of the
proteins, and the terms $1-c_j, \ j=1,2,3,$
indicate the inactive forms (after
non-dimensionalization and assuming that the total concentration of
each of the proteins is $1$).
Here the parameters $ a_1, a_2, a_3,
   b_1, d_1, e_1,
   b_2, d_2, e_2,
   b_3,
d_3, e_3,   \mu, k $
are strictly positive.

Let
$\Omega:= \{c\in \mathbb{R}^3: 0 < c_j < 1 , \ j=1,2,3 \}$
    denote the open unit cube,
    the domain where this model holds.
Now we show how to  apply Theorem  \ref{thm:intro}
    on $\Omega$ to conclude that:
\begin{quote}
{\it
There is a unique equilibrium in $\Omega$ for any
    choice of the strictly positive  parameters,
$ a_1, a_2, a_3,
   b_1, d_1, e_1,
   b_2, d_2, e_2,
   b_3,
d_3, e_3,   \mu, k $.}
\end{quote}

\noindent
    {\bf Proof:} \
First the Jacobian has the form (\ref{eq:Amatrix}).
Thus (\ref{eq:detAS})
implies that the determinant is strictly negative for all strictly
positive parameters and
concentrations $c\in \Omega$.
The proof follows the same outline  as Example \ref{ex:AS1}.
Now we check the required items:\\ \\
{\it No equilibria on the  boundary
of the unit cube:} \
   Suppose there
is an equilibrium $c$ on the boundary of $\Omega$.
Then the equilibrium equations imply that
\begin{enumerate}
\item
If $c_1=0$ then (\ref{amapk1}) forces $c_1=1$. Contradiction.

\item
If $c_1=1$ then (\ref{amapk1}) forces $c_1=0$. Contradiction.

\item
If $c_2=0$ then (\ref{amapk2}) forces $c_1=0$. Contradiction as above.

\item
If $c_2=1$ then  (\ref{amapk2}) forces $c_2=0$. Contradiction.

\item
If $c_3=0$ then (\ref{amapk3}) forces $c_2=0$.
    Contradiction as above.

\item
If $c_3=1$ then (\ref{amapk3}) forces $c_3=0$. Contradiction.
\end{enumerate}

\noindent
{\it Initializing:} \
\cite{ASprept}
  proves that there are choices of parameters compatible with
this model for which
 there is a unique stable equilibrium point in $\Omega$.
Alternatively,  one can  compute  for a simple
choice of parameters that there is  a unique
  positive equilibrium.

The discussion exactly as before implies that
there is a unique equilibrium for all 
fixed strictly positive values of the parameters.
\qed
\end{exam}

We mention here that the question of how one finds
good intializations for the rate constants might be
a topic for further research.
The goal would be to find methods for
systematically selecting rate constants
that produce systems whose equilibria can be determined
by analytic means.
We have not explored this topic at all.

\subsection{The theory of Arcak and Sontag}
\label{sec:ASsummary}

Now we shall make some general comments on
\cite{arcson06, ASprept}.
There were four chemical
   reaction examples
presented  in the two papers \cite{arcson06, ASprept}.
So far we have treated two of the four here in this section.
The third example, Example 2 in \S4 of \cite{ASprept}, is a small variant of
Example \ref{mapk} above  and it  can be treated in a similar manner to
that example. In particular, it  has a Jacobian of the form
(\ref{eq:Amatrix}). 
We now  turn to the fourth example of Arcak and Sontag.

\begin{exam}
\label{exAS3} This is  Example 3 in \S4 of \cite{ASprept} which we do not
describe in detail, since it requires about a page.  While its
Jacobian
  does not have the form (\ref{eq:Amatrix}), it is easy to
analyze
  (using Mathematica) and what we found is
     that the  determinant of the Jacobian
of $f$ is positive at all  strictly positive $c $.
Thus the theory described here applies provided suitable boundary
behavior holds.
Boundary behavior was not possible to determine
since the example was a rather general class
whose boundary behavior was not specified.
In a particular case where more information is specified
one might expect that this could be done.
\qed
\end{exam}

Arcak and Sontag \cite{ASprept} present
a general theory which contains  the  examples considered in this section
and which does not match up simply with ours.
Their theory assumes an equilibrium exists (we do not).
It places global restrictions on the equilibrium which guarantee
that it is a unique globally stable equilibrium (we address uniqueness
but not stability).
However, while the Arcak-Sontag theory is different than
ours, we do point out in this section
that all four of their chemical examples
have Jacobians whose determinant sign does not depend on chemical
concentration, so our approach applies directly,
and with a bit of attention to boundary behavior, gives
existence of a unique positive equilibrium.
However, we do not obtain the very impressive global
stability in \cite{ASprept}.

\section{Degree and homotopy of maps}
\label{sec:degree}
The proof of Theorem \ref{thm:intro} and other results in this
paper
is based on classical degree theory.
The degree of a function is invariant if we continuously
deform
(homotopy) the function and we use that to advantage in this paper.

Now we give the setup.  If $\Omega \subset {\mathbb R}^n$
  is a
bounded domain, and if a smooth
(once continuously differentiable) function
$f:\overline\Omega \to {\mathbb R}^n$
has no degenerate zeros, and has no zeros on the boundary
of $\Omega$, then the {\it topological degree with respect to zero
of $f$} (or simply the \emph{degree of $f$}) equals
\beq\label{def_degree}
\deg(f) = \deg(f,\Omega) = \sum\limits_{c \in Z_f} \text{ sgn }
\left(\det \left( \frac{\partial f}{\partial c} (c)\right) \right),
\eeq
where $\text{sgn}:{\mathbb R} \to \{-1,0,1\}$ is the sign function,
    $Z_f$ is the set of zeros of $f$ in $\Omega$, and $c^\ast$ is
  a
degenerate point means $\det\left(\frac{\partial f}{\partial c}(c^\ast)\right)=0$.
The degree of a map naturally extends from nondegenerate
   smooth functions to
 continuous functions $f:\overline\Omega \to {\mathbb R}^n$.
For this, one can approximate $f$ uniformly with smooth functions $F_k$
that have no degenerate zeros and no zeros on the boundary
of $\Omega$, and then define the degree of $f$ as the limit
of the degrees of $F_k$. The key fact is:
this construction of the degree of $f$ is independent of
the approximates $F_k$.  Fortunately, we shall only need to compute
deg$(f)$ on smooth nondegenerate $f$.
For a quick account of this theorem and the main properties of degree,
see Ch 1.6A of \cite{Berger}.

\emph{Homotopy invariance of the degree}
is the following well known property:

\begin{theorem}\label{thm:degree}
Consider some bounded domain $\Omega \subset {\mathbb R}^n$
    and a continuously varying family of smooth functions
    $f_\lambda : \overline\Omega \to {\mathbb R}^n$ for $\lambda\in [0,1]$,
    such that $f_\lambda$ does not have any zeros on the boundary
    of $\Omega$ for all $\lambda\in [0,1]$.
    Then $\deg(f_\lambda)$ is constant for all $\lambda \in [0,1]$.
\end{theorem}

Now we give a slightly more general theorem than
  Theorem
\ref{thm:intro} stated in the introduction.

\begin{theorem}
\label{thm:degreeDet}
Suppose $\Omega$ and $f_\lambda, \ \lambda \in [0,1]$, are
as in Theorem \ref{thm:degree}.
Then for any $\lambda \in [0,1]$ such that
$\det \left( \frac{\partial f_\lambda}{\partial c} (c) \right) \neq 0$
for all $c \in \Omega$,
the number of zeros of $f_\lambda $ in $  \Omega$ must equal
the absolute value of the degree of $f_\lambda $ in $ \Omega$, which
equals the absolute value of the degree of $f_{\lambda'}$ for any
$\lambda'\in [0,1]$.
\end{theorem}

\noindent
{\bf Proof:}
If $\lambda \in [0,1]$ is such that the determinant
$\det({\partial f_\lambda}/{\partial c})$ does not vanish
in $\Omega $, then
      $\text{sgn} \left( \det({\partial f_\lambda}/{\partial c})\right)$ is
independent of $c$.
This implies that the zeros of $f_\lambda$ are nondegenerate and, by
the formula for the degree of $f_\lambda$, that
$ |\deg (f_\lambda)|$ equals the number of zeros
     of $f_\lambda$ in $\Omega$. The fact that the degree
     does not vary with $\lambda$ is immediate from Theorem \ref{thm:degree}.
\qed


\begin{rem}
For $|\deg(f_\lambda)|$ to count the number of zeros of $f_\lambda$
in $\Omega$,
$\text{sgn}\left( \det \left( \frac{\partial f_\lambda}{\partial c}
  (c^*) \right)\right)$ need only be the same
for  all zeros
  $c^*$ in $\Omega$, not for all $c \in \Omega$.
Sadly this weakening of hypotheses
  is hard to take advantage of in
practice.
\end{rem}

\begin{rem}
  From the viewpoint of numerical calculation, Theorem
\ref{thm:degreeDet} strongly suggests that if the
no boundary zeros hypothesis holds,  and (DetSign) holds
for $f=f_\lambda$ at one value of $\lambda=\lambda_1$,
and if one can calculate all zeros of $f_\lambda$
at some other value of $\lambda= \lambda_2$ where (DetSign) also holds,
then we can determine the number
of zeros at $\lambda=\lambda_1$.
Indeed, often we can find a $\lambda_2$ for which $f_{\lambda_2}$
is ``simple" in the sense that all zeros for $\lambda_2$ are
non-degenerate  and the zeros can be readily computed along  with the
Jacobians there,
 and consequently
$\deg(f_{\lambda_2})$
can be computed.
The import for numerical calculation is that
finding a single equilibrium is often not so
onerous.
After finding one equilibrium one typically makes
a new initial guess and tries to find another. Knowing if one
has found all of the equilibria is the truly daunting task, since
it is nearly impossible to ensure this by experiment.
Thus theoretical results (hopefully those here) help with
this very difficult computational question.
\end{rem}

\section{Mass-dissipating dynamical systems}
\label{sec:mass}

In this section we consider a general
dynamical system model which includes
the more specific dynamics of
conservative chemical reaction networks, augmented with
inflows and outflows, as  described in the next section.
In chemical engineering, the latter is commonly refered
to as dynamics that goes with a continuous flow stirred tank
reactor (CFSTR).
In biochemistry, one may view this as a model for intracellular
behavior with production and degradation, or with inflow and outflow
across the cell boundary.
Here all species components are subject to inflow and outflow,
however, to approximate the conservation of some
species such as enzymes, one may take the associated inflow rate
value in $c_{in}$ and degradation factor in $\lout$ to be arbitrarily small,
if desired.

In preparation for defining a dynamical system on the orthant
$\bRplusc$,
we consider a smooth  function
$g :\bRplusc  \to {\mathbb R}^n$,
where $g$ has the property that
for each $j \in \{ 1,...,n \}$,
the $j^{\rm th}$ coordinate of $g(c)$ is nonnegative
whenever
     the $j^{\rm th}$ coordinate of
    $c \in  {\mathbb R}^n_{\ge 0}$ is zero.
Consider the  dynamical system associated with this function
given by
\begin{equation} \label{rbasic}
\dot{c}= g(c) \quad \hbox{ for }  c\in \Rnn.
\end{equation}
This
dynamical system (\ref{rbasic})
is called {\it positive-invariant}
because of the condition on $g$.
This guarantees that the dynamics leaves the orthant
$\bRplusc$ invariant.
Given $m\in \mathbb{R}_{>0}^n$, the dynamical system (\ref{rbasic})
is called {\it mass-dissipating with respect to $m$}
if 
\beq
\label{eq:massDiss}
m \cdot g(c) \le 0
\eeq
  for all $c \in \bRplusc$;
it is called {\it mass-dissipating} if it is mass-dissipating
with respect to $m$ for some $m \in {\mathbb R}^n_{>0}$.
In this case, on $\Rnn$,
\begin{equation}
\frac{d(m\cdot c) }{dt}= m\cdot g(c) \leq 0.
\end{equation}

Now we consider
the dynamical system (\ref{rbasic})  augmented
with inflows and outflows:
\be \label{dissipating}
\dot{c} = c_{in} -\lout c + g(c),
\ee
where
$\lout$ is  an $n\times n$ diagonal matrix with strictly positive entries
on the diagonal.  We interpret
the term $c_{in} \in {\mathbb R}^n_{>0}$ as a constant {\it inflow rate},
and the term $\lout c$ as an {\it outflow rate} which for each component
is proportional to the concentration of that component.
It is easy to check that with this augmentation, the dynamics still
leaves the orthant $\Rnn$ invariant.
However, the mass-dissipating property is only inherited at large
values of the concentration $c$.

We are now ready to state our main theorem in this context.

\begin{theorem}
\label{thm:dissip}
Let $c_{in}, m \in {\mathbb R}^n_{>0}$,  and $\lout$ be an
$n\times n$ diagonal matrix with strictly positive diagonal entries.
Consider a smooth function
$g : { \bRplusc} \to {\mathbb R}^n$
such that the dynamical system (\ref{rbasic})
is positive-invariant and mass-dissipating with respect to $m$.
Define
$$f(c): = c_{in} -\lout c + g(c) \ \ \hbox{ for } c \in \bRplusc. $$
Then the augmented system (\ref{dissipating}),  with inflows
and outflows,
has no equilibria on the boundary of
$\Rnn$, and if
$\det\left(\frac{\partial f}{\partial c}\right)\not= 0$ on
$\bRplus $,
then there is  exactly one equilibrium point for this system in $\Rnp$.
\end{theorem}


\begin{proof}
It suffices to prove that $f$ has no zeros on the boundary of $\Rnn$
and if
$\det\left(\frac{\partial f}{\partial c}\right)\not= 0$ on
$\bRplus$, then $f$ has exactly one zero in $\Rnp$.

Define
$$f_\lambda(c):=  c_{in} - \lout c + \lambda g(c),\ \ \hbox{for } c\in \Rnn,\ \lambda\in [0,1].$$
Fix $M> m\cdot c_{in}$ and let
$${\Omega}_M=\{ c\in \bRplus : m\cdot(\lout c) < M\}. $$
Then ${\Omega}_M$ is a bounded domain and
$\{ f_\lambda: \lambda \in [0,1]\}$ is
a continuously varying family of smooth
functions on $\overline{\Omega}_M$. For $j= 1, \ldots, n,$
consider $c^j \in \overline{\Omega}_M$
such that the $j^{\rm th}$ coordinate of $c^j$ is zero.
Then the $j^{\rm th}$ coordinate of $f_\lambda(c^j)$
must be strictly positive, because the $j^{\rm th}$
coordinate of $c_{in}$ is strictly positive, and the $j^{\rm th}$
coordinate of $g(c^j)$ is nonnegative,
by the positive-invariance assumption.
Therefore $f_\lambda$ has no zeros on
the sides of ${\Omega}_M$, i.e., on  $\overline{\Omega}_M \cap
\partial\bRplusc$.
Also, we have
\beq
\label{eq:storage}
  m \cdot f_\lambda(c)  =  m \cdot c_{in} -
m \cdot (\lout  c ) +\lambda m \cdot g(c)
    \le  m \cdot  c_{in} - m \cdot(\lout  c) < 0
    \eeq
for all $c\in \overline{\Omega}_M$ such that $ m \cdot(\lout c)  = M$.
Here we have used the mass-dissipating property of $m$ for the first
inequality and the fact that $M> m\cdot c_{in} $ for the second inequality.
It follows that 
$f_\lambda$ has no zeros on
the outer boundary of ${\Omega}_M$, i.e., on $ \{ c \in \bRplus : m
\cdot (\lout c ) = M\}$.
Thus, $f_\lambda$ has no zeros on the boundary of ${\Omega}_M$
for all $\lambda \in [0,1]$. Then, by Theorem \ref{thm:degree},
the degree of $f_\lambda $ on ${\Omega}_M$,
$\deg(f_\lambda,{\Omega}_M)$, is constant for all $\lambda \in [0,1]$.
Next we observe that $c^*={(\lout)}^{-1}c_{in}$ is the unique zero of
$f_0$ and is inside ${\Omega}_M$, and $\frac{\partial f}{\partial c}= -\lout$, and so by (\ref{def_degree}), we obtain
$\deg(f_0, {\Omega}_M)= \text{sgn}(\det(-\lout))=(-1)^{n}$.
Hence, by Theorem \ref{thm:degreeDet},
if $\det\left(\frac{\partial f}{\partial c}\right)
=  \det \left(\frac{\partial f_1}{\partial c}\right)\not= 0$
on ${\Omega}_M$, then $f=f_1$
has exactly one zero in ${\Omega}_M$.

Since $M>m\cdot c_{in}$ was arbitrary and the sets
$\overline {\Omega}_M:M>m\cdot c_{in}$  fill out $\bRplusc$,
it follows that
$f$ has no zeros on the boundary of $\bRplusc$.  If furthermore,
$\det\left(\frac{\partial f}{\partial c}\right)\not= 0 $ on all of
$\bRplus$, then it follows that $f$ has exactly one zero in $\bRplus$.
\end{proof}

\begin{rem}
\label{rem:massAtoms}
A special case of mass-dissipating is \it mass-conserving, \rm
namely
$ m\cdot g(c)=0$.
For dynamical systems associated
to chemical reaction networks  this has a natural interpretation.
Indeed, the dynamics of
chemical concentrations resulting from chemical interactions
among several types of
molecules will be mass-conserving whenever there exists a
    mass assignment for each
chemical species which is conserved by each reaction,
or whenever each chemical
species (or molecule) is made up of atoms that are also conserved
by each reaction.
More generally, the dynamics will be mass-dissipating
whenever no reaction produces
more mass than it consumes,
respectively,  produces more atoms than it consumes.
Mass conservation implies $\det (\frac{ \pt g}{\pt c}) =0$,
since $m \cdot \frac{ \pt g}{\pt c} = 0$ when $m\cdot g =0$.
Thus augmenting with outflows is required to make the
hypothesis on the sign of $ \det (\frac{ \pt f}{\pt c})$
in our theorems meaningful.
The paper  \cite{HKG08}, which builds on the current one,  introduces
a more general determinant that applies when there are no outflows
(or only some outflows). This then helps one  count equilibria  in a manner
generalizing what we have done here.
\end{rem}

\begin{rem}
Theorem \ref{thm:dissip} still holds with a much
less restrictive definition of mass-dissipating,
e.g., by replacing ``$m$" with the gradient ``$\nabla L$" for
an appropriate class of functions $L:\Rnn\to \mathbb{R}_{\geq 0}$.
Here mass or atom count is behaving like what is called
{\it storage function} in engineering systems theory,
see \cite{K02}.
Indeed the inequality
 $m \cdot f_\lambda(c) \le  m \cdot  c_{in} - m \cdot (\lout  c)$
derived in (\ref{eq:storage}) is what is called a dissipation
inequality on the storage function $c\to m\cdot c$, which in fact
also plays the role of a  ``running cost".
\end{rem}

\section{Dynamics of chemical reaction networks}
\label{sec:dynam}

We now introduce the standard terminology of
Chemical Reaction Network Theory
(see \cite{HornJackson72, Fein72, Fein95, CF05}).
A chemical reaction network is usually specified by a
finite set of reactions $\mathscr R$ involving a finite set of chemical species
$\mathscr S$.

For example, a chemical reaction network with
two chemical species $A_{1}$ and $A_{2}$ is schematically
given in the diagram

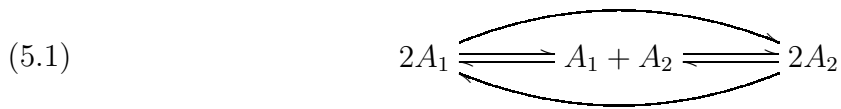
\begin{equation}\label{ex1}
\xymatrix @ +1pc{
2A_{1} \ar@<.3 ex>@^{->}[r]\ar@/^/ @/^1.5pc/ @{-^{>}}[rr]
&A_{1}+A_{2}\ar @<.3 ex>@^{->}[l]\ar @<.3 ex>@^{->}[r]
&2A_{2}\ar  @{-^{>}}@/^1.5pc/[ll]\ar@<.3 ex>@^{->}[l]}
\end{equation}

The dynamics of the state of this chemical system is defined
in terms of functions $c_{A_{1}}(t)$ and $c_{A_{2}}(t)$
which represent the  concentrations of the species
$A_1$ and $A_2$ at time $t$.
The occurrence of a chemical reaction causes changes in concentrations;
for instance, whenever the reaction $A_{1}+A_{2}\to 2A_{1}$ occurs,
the net gain is a molecule of $A_{1}$,
whereas one molecule of $A_{2}$ is lost.
Similarly, the reaction $2A_{2}\to 2A_{1}$ results in
    the creation of two molecules of $A_{1}$ and the loss
    of two molecules of $A_{2}$.

A common assumption is that the rate of change of the concentration
of each species is governed by {\it mass-action kinetics}
\cite{HornJackson72, Fein72, Fein79, FeinZero, Fein95, Sontag,
CF05, CF06, CTF06}, i.e., that each reaction takes place at a rate
that is proportional to the product of the concentrations of the
species
    being consumed in that reaction. For example, under the
mass-action kinetics assumption, the contribution of the reaction
$A_{1}+A_{2}\to 2A_{1}$ to the rate of change of $c_{A_{1}}$ has
the form $k_{A_{1}+A_{2}\to 2A_{1}}c_{A_{1}}c_{A_{2}}$,
    where
$k_{A_{1}+A_{2}\to 2A_{1}}$ is a positive number called the
    {\it
reaction rate constant}. In the same way, the reaction
    $2A_{2}\to
2A_{1}$ contributes the negative value
    $-2k_{2A_{2}\to
2A_{1}}c_{A_{2}}^{2}$ to the rate of change of
    $c_{A_2}$.
Collecting the contributions of all the reactions,
    we obtain the
following dynamical system associated
    to the
chemical reaction network depicted in (5.1):

\begin{eqnarray}\label{diff_eq}
&\dot c_{A_{1}}
    =&-k_{2A_{1}\to A_{1} +
A_{2}}c_{A_{1}}^{2}+k_{A_{1}
    + A_{2}\to 2A_{1}}c_{A_{1}}c_{A_{2}}
- k_{A_{1}+A_{2}
    \to 2A_{2}}c_{A_{1}}c_{A_{2}}\\ \nonumber
&
    &+\, k_{2A_{2}\to A_{1}+A_{2}}c_{A_{2}}^{2}
    -
2k_{2A_{1}\to 2A_{2}}c_{A_{1}}^{2}+2k_{2A_{2}
    \to
2A_{1}}c_{A_{2}}^{2} \\\nonumber
&\dot c_{A_{2}}
    =&k_{2A_{1}\to A_{1} + A_{2}}c_{A_{1}}^{2}
    -
k_{A_{1} + A_{2}\to 2A_{1}}c_{A_{1}}c_{A_{2}}
    + k_{A_{1}+A_{2}\to
2A_{2}}c_{A_{1}}c_{A_{2}}\\ \nonumber
&
    &-\, k_{2A_{2}\to A_{1}+A_{2}}c_{A_{2}}^{2}
    +
2k_{2A_{1}\to 2A_{2}}c_{A_{1}}^{2} - 2k_{2A_{2}
    \to
2A_{1}}c_{A_{2}}^{2} \\\nonumber
\end{eqnarray}

The objects on both sides of the reaction arrows
(i.e., $2A_{1}, A_{1}+A_{2}$, and $2A_{2}$) are called
    \emph{complexes} of the reaction network.
    Note that the complexes are non-negative integer combinations
    of the species. On the other hand, we will see later that it
    is very useful to think of the complexes as (column) vectors in $\RR^n$,
     where $n$ is the number of elements of $\mathscr S$,
    via an identification of the set of species with
the standard basis of
    $\RR^n$, given by a fixed ordering of the
species.
    For example, via this identification, the complexes above
become
    $2A_{1}=\begin{bmatrix} 2 \\ 0 \\\end{bmatrix}$,
$A_{1}+A_{2}=\begin{bmatrix} 1 \\ 1 \\\end{bmatrix}$, and
$2A_{2}=\begin{bmatrix} 0 \\ 2 \\\end{bmatrix}$.
    We can now
formulate
a general setup which includes
many situations, certainly those above.

\subsection{The general setup}

Now we present basic definitions and illustrate them.

\begin{definition}
A \textit{chemical reaction network} is a
triple $(\mathscr{S}, \mathscr{C}, \mathscr{R})$,
where $\mathscr{S}$ is a set of $n$ chemical \textit{species},
$\mathscr{C} $ is a finite set  of vectors in $\bRplusc$ with
nonnegative integer entries  called
the set of \textit{complexes},
and
%
%
$\mathscr{R} \subset \mathscr{C} \times \mathscr{C}$
is a finite set of relations
between elements of $\mathscr{C}$,
denoted $y\to y'$ which  represents the set of
    \textit{reactions} in the network.
Moreover,
the set    $\mathscr{R}$ cannot contain
elements of the form $y\to y$;
    for any $y\in \mathscr{C}$ there
exists some $y'\in\mathscr{C}$
    such that either $y\to y'$ or
$y'\to y$; and the union of the
supports of all $y\in\mathscr{C}$ is $\mathscr{S}$,
where the \textit{support} of an element $\alpha\in \mathbb R^n$ is
$\supp(\alpha)=\{ j : \alpha_{j}\neq 0\}$.
To each reaction $y\to y'\in \mathscr{R}$, we associate a reaction
vector given by $y'-y$.
\end{definition}

The last two constraints of the definition amount to
requiring that each complex appears in at least
one reaction, and each species appears in at least one complex.
For the system (\ref{ex1}), the set of species is
$\mathscr{S}=\{A_{1}, A_{2}\}$, the set of complexes is
$\mathscr{C}=\{2A_{1}, A_{1}+A_{2}, 2A_{2} \}$ and the set of
reactions is $\mathscr{R}=\{2A_{1}
    \rightleftharpoons A_{1}+A_{2},
A_{1}+A_{2}
    \rightleftharpoons 2A_{2}, 2A_{2}
    \rightleftharpoons
2A_{1}\}$,
    and consists  of 6 reactions, represented as three
reversible reactions.


In examples we will often refer to a chemical reaction network by
specifying $\mathscr{R}$ only, since $\mathscr{R}$ encompasses
all of the information about the network.
In the sequel we shall sometimes simply say reaction network in place of chemical reaction network.

\begin{definition}
A \emph{kinetics} for a reaction network
    $(\mathscr{S},
\mathscr{C}, \mathscr{R})$ is an assignment to
    each reaction $y\to
y' \in \mathscr{R}$ of a \emph{reaction rate function}
$$
K_{y\to y'} : {\RR}_{\ge 0}^n  \to {\RR}_{\ge 0}^n.
$$
By a \emph{kinetic system}, which we denote by
    $(\mathscr{S},
\mathscr{C}, \mathscr{R}, K)$, we mean a reaction
    network taken
together with a kinetics.
\end{definition}

For each \emph{concentration} $c \in {\RR}_{\ge 0}^n$,
the nonnegative number $K_{y\to y'}(c)$
is interpreted as the occurrence rate
of the reaction $y \to y'$ when the chemical mixture
has concentration $c$.
Hereafter, we suppose that reaction rate functions
are smooth on ${\RR}_{\ge 0}^n$,
and that $K_{y\to y'}(c) = 0$
whenever $\supp(y) \not\subset \supp(c)$.
Although it will not be important in this article,
it is natural to also require that,
for each $y \to y' \in \mathscr{R}$
the function $K_{y\to y'}$ is strictly positive
precisely when $\supp(y) \subset \supp(c)$,
i.e., precisely when the concentration $c$
contains at non-zero concentrations
those species that appear in the reactant complex $y$.

\begin{definition}
The \emph{species formation rate function} for a kinetic system
$(\mathscr{S}, \mathscr{C}, \mathscr{R}, K)$
is defined by
$
r: {\RR}_{\ge 0}^n  \to {\mathbb R}^n $
where
$$ r(c)
    = \sum_{y\to
y'\in\mathscr{R}}K_{y\to y'}(c)(y'-y) \ \ \hbox{for } \ c\in \Rnn .
$$
The
{\emph associated dynamical system}    for the
kinetic system $(\mathscr{S}, \mathscr{C}, \mathscr{R}, K)$ is
\begin{equation}\label{diff_chem}
       \dot c = r(c) = \sum_{y\to y'\in\mathscr{R}}K_{y\to y'}(c)(y'-y),
\end{equation}
where $c\in \bRplusc$ is the nonnegative vector of species concentrations.
\end{definition}

The interpretation of $r(\cdot)$ is as follows:
if the chemical concentration is $c \in {\RR}_{\ge 0}^n$,
then $r_j(c)$ is the production rate of species $j$
due to the occurrence of all chemical reactions.
To see this, note that
$$
r_j(c) = \sum_{y\to y'\in\mathscr{R}}K_{y\to y'}(c)(y'_j-y_j),
$$
and that $y'_j-y_j$ is the net number of molecules of species $j$
produced with each occurrence of reaction $y\to y'$.
Thus, the right hand side of the equation above is the sum of
all reaction occurrence rates, each weighted by the net gain
in molecules of species $j$ with each occurrence of the
corresponding reaction.
Note that $r_j(c)$ could be less than zero, in which case
$|r_j(c)|$ represents the overall rate of consumption of species $j$.

\subsubsection{Special case: Mass-action kinetics}

\begin{definition}
A \textit{mass-action system} is a quadruple
    $(\mathscr{S},
\mathscr{C}, \mathscr{R}, k),$
    where $(\mathscr{S}, \mathscr{C},
\mathscr{R})$ is a
chemical reaction network and
$k = (k_{y\to y'})$
is a vector of \textit{reaction rate constants},
so that the reaction rate function
$
K_{y\to y'} : {\RR}_{\ge 0}^n  \to {\RR}_{\ge 0}^n,
$
for each reaction ${y\to y'}\in \mathscr R,$ is given by
\emph{mass-action kinetics}:
$$
K_{y\to y'}(c)    = k_{y\to y'}c^y \qquad
where \ \ \ c^y = \displaystyle \prod_{i =1}^n c_i^{y_i}.
$$
(Here we adopt the convention that $0^0 =1$.)
The
\emph{associated mass-action dynamical system} is
\begin{equation}\label{mass_act}
       \dot c= \sum_{y\to y'\in\mathscr{R}}k_{y\to y'} c^y (y'-y).
\end{equation}
\end{definition}
In the vector equation (\ref{mass_act}),
the total rate of change of the vector of concentrations $c$
is computed by summing the contributions of all the reactions
in $\mathscr R$. Each reaction $y\to y'$ contributes proportionally
to the product of the concentrations of the species in
its {\it source} $y$, that is, $c^{y}$, and also proportional to the
number of molecules gained or lost in this reaction.
   Finally, the
proportionality factor is $k_{y\to y'}.$
   For example, we can
rewrite (\ref{diff_eq}) in the vector form
   (\ref{mass_act}) as
\begin{alignat}{4}
\begin{bmatrix}
\dot c_{1}\\
\dot c_{2}\\
\end{bmatrix}&=&\ k&_{2A_{1}\to A_{1}+A_{2}}c_1^{2}
\begin{bmatrix}
-1\\
1\\
\end{bmatrix}+k_{A_{1}+A_{2}\to 2A_{1}}c_1c_2
\begin{bmatrix}
1\\
-1\\
\end{bmatrix}+k_{A_{1}+A_{2}\to 2A_{2}}c_1c_2
\begin{bmatrix}
-1\\
1\\
\end{bmatrix}\\\nonumber
& &+k&_{2A_{2}\to A_{1}+A_{2}}c_2^{2}
\begin{bmatrix}
1\\
-1\\
\end{bmatrix}+k_{2A_{1}\to 2A_{2}}c_1^{2}
\begin{bmatrix}
-2\\
2\\
\end{bmatrix}+k_{2A_{2}\to 2A_{1}}c_2^{2}
\begin{bmatrix}
2\\
-2\\
\end{bmatrix}.
\end{alignat}

\subsection{Mass conservation}

Now we see in terms of this setup
how one obtains mass conservation as defined
in \S \ref{sec:mass}.

\begin{definition} The \emph{stoichiometric subspace}
$S \subset \RR^n$ of a reaction
network $(\mathscr{S}, \mathscr{C}, \mathscr{R})$
is the linear subspace of $\RR^n$
spanned by the {\it reaction vectors} $y'-y$,
for all reactions $y\to y' \in \mathscr{R}$.
\end{definition}

Note that, according to (\ref{diff_chem}), for a given
value of $c$, the vector $\dot c$
   is a
linear combination of the reaction vectors.
   This implies that each
\emph{stoichiometric compatibility  class}
   $(c_0 + S) \cap
\RR_{\ge0}^n$ is an invariant set for
   the dynamical system
(\ref{diff_chem}) with initial condition $c_0 \in \Rnn$.

\begin{definition} A reaction network
$(\mathscr{S}, \mathscr{C}, \mathscr{R})$
is called \emph{conservative} if there exists some
positive vector $m \in \mathbb{R}_{>0}^n$
which is orthogonal to all its reaction vectors, i.e.,
$$ m\cdot (y'-y)=0 $$
for all reactions $y\to y'$ in $\mathscr{R}$.
Then $m$ is called a {\it conserved mass vector}.
\end{definition}

\begin{rem}
Each trajectory of a conservative
reaction network is bounded.
A conservative reaction network can have no inflows
or outflows (see the definition of inflow and outflow in the
next section).
\end{rem}

\subsection{Main results}

We  now consider conservative reaction networks augmented
with inflow and outflow reactions (for each of the species).
An inflow reaction is a reaction of the form $0\to A$
and an outflow reaction is one of the form $A\to 0$, where
$A$ is a species.
The reaction vector $y'-y$ associated with an inflow
reaction for species $j$ is the vector containing all zeros, except
that it has a one in the $j^{th}$ position.
The reaction vector  associated with an outflow
reaction for species $j$ is the negative of  the vector associated with
an inflow reaction for that species.
Here, for the kinetics associated with the inflows and outflows,
we assume that the reaction rate function for each
inflow reaction is a positive constant and the value of the reaction rate
function for each outflow reaction is a positive constant
times the concentration of the species flowing out.
The latter corresponds to degradation of each species at a rate
proportional to its concentration.
The following theorem may be used to determine the number
of equilibria for conservative reaction networks augmented by
such inflows
and outflows.  It requires a positive determinant condition
and is
the analog of Theorem 4.1 in this context.


\begin{theorem}
\label{thm:mm}
Consider some conservative reaction network
$(\mathscr{S}, \mathscr{C}, \mathscr{R})$ with conserved mass vector $m$.
Let $K$ be a kinetics for this network with associated species formation
rate function $g$.
Consider an augmented kinetic system
$(\mathscr{S},\tilde{\mathscr{C}},\tilde{\mathscr{R}},\tilde{\mathscr{K}})$
obtained by adding
inflow and outflow reactions for all species
so that the associated dynamical system is:
\beq
\label{eq:augK}
\dot c = r(c):= c_{in} -\lout c + g(c),
\eeq
where $c_{in} \in \bRplus$ and $\lout $ is an $n\times n$ diagonal matrix
with strictly positive diagonal entries.
Suppose that
$$
\det \left( \frac{\partial r}{\partial c} (c) \right) \neq 0, \ \
$$
for all
    $c \in {\mathbb R}^n_{> 0}
$.
Then  the dynamical system (\ref{eq:augK})
has exactly one equilibrium
$c^*$ in $\bRplus$ and no equilibria on the boundary of $\bRplusc$.
\end{theorem}

\begin{proof}
We want to apply Theorem \ref{thm:dissip}.
The function $g $ is given by the right member of (\ref{diff_chem}),
where the functions $K_{y\to y'} $ are all smooth and
have the property that $K_{y\to y'}(c)=0 $ whenever
$\supp(y)\not\subset \supp(c)$.
Consequently, $g$ is smooth and, whenever
$c\in \bRplusc$ is such that $c_j=0$ for some $j$,
then we have
$$g_j(c) \geq \sum_{y\to y' \in \mathscr{R}} K_{y\to y'}(c)(-y_j)=0,$$
because  $K_{y\to y'}(c)=0$ whenever $y_j>0$ and $c_j=0$,
by the support property of $K$.
It follows that the dynamical system (\ref{eq:augK})
is positive-invariant.
Furthermore, the system is mass-dissipating, since
$$m\cdot g(c) =
\sum_{y\to y' \in \mathscr{R}} K_{y\to y'}(c)\, m\cdot(y-y') =0,$$
by the assumption that
the reaction network $(\mathscr{S}, \mathscr{C}, \mathscr{R})$
is conservative.
The conclusion then follows immediately from Theorem
\ref{thm:dissip}.
\end{proof}

\begin{rem}
The results described above use the assumption that all species
have inflows, in order to conclude that there are no equilibria
on the boundary of $\Rnp$. On the other hand, for very large classes
of chemical reaction networks described in \cite{ADSprept},
this assumption is actually not needed in order to rule out
the existence of such boundary equilibria (an observation
by David Anderson University of Wisconsin \cite{DavidAnderson}).
\end{rem}

\begin{rem}
The paper \cite{BDB07}, for the case of ``nonautocatalytic reactions",
gives a condition on the ``stoichiometric matrix"
(in our terminology
the matrix whose columns are the vectors $y'-y$ for $y \to y' \in \mathscr R$)
  which is necessary and sufficient for  the determinant of the Jacobian of $r$
to be of one sign for all concentrations and all $K_{y \to y'}$
which are monotone increasing in each variable.
Our theory is less restrictive as illustrated by
Example \ref{ex:CF1}.
\end{rem}

\section{Applications}
\label{sec:appl}

In practice, most dynamical system models of biochemical reaction
networks contain a large number of unknown parameters.
These parameters correspond
to reaction rates and other chemical
properties of the reacting species.
In this section
we treat a variety of examples of such models and illustrate
the use of Theorems \ref{thm:dissip} and \ref{thm:mm}.
In some of these examples, the determinant of the Jacobian
$\det\left(\frac{\pt r}{\pt c}\right)$ is of one sign
everywhere on the open orthant for all parameters
and in some it is not.
Even in the latter cases we describe  ways to find
classes of rate functions
for which there exists a unique positive equilibrium.

The first subsection assumes mass-action kinetics
and defines (reminds)
the reader of the Craciun-Feinberg ``determinant expansion"
via an example.
Critical is the sign of each term in the expansion
 and whether
all terms have the same sign or  miss this by ``a little",
namely, only one or two terms in the determinant expansion
has a  coefficient with an anomalous sign.
Here we point out that all examples
in \cite{CF05, CF06, CTF06}  have
at most one or two anomalous signs.
When there are no anomalous signs, these papers
show that any positive equilibrium is unique for all
parameter values, and they  
develop and use graph theoretical methods 
for determining when there are no anomalous signs.
Here, for cases of few anomalous signs,
we propose and illustrate a technique for identifying
parameter values for which a positive equilibrium
exists and is unique.
The paper, \cite{HKG08}, subsequent to this one, 
gives ways of counting the number of anomalous signs
in terms of graphs associated to
a chemical reaction network.

In the second subsection,
we continue with the general framework of \S \ref{sec:dynam}, and move beyond
mass-action kinetics to rate functions satisfying certain monotonicity
conditions 
(see Definitions \ref{def:incrCI} and \ref{def:incrSM}).
The  weaker condition,
Definition \ref{def:incrSM}, holds for many biochemical reactions
and allows one to make sense of the signs which
occur in the  determinant expansion.
Hence one can  apply the methods in this paper.

The number of anomalous signs can be determined for the examples in this
section
using: (a) the graph-theoretic methods  of Craciun and Feinberg \cite{CF05, CF06} 
when there are no anomalous signs and the kinetics are of mass-action type,
and (b)
symbolic computation of the determinant of the 
Jacobian using Mathematica when there are some anomalous signs 
or the  kinetics are general (not necessarily mass-action). 
The reader will find Mathematica
notebooks at 
$$\hbox{http://www.math.ucsd.edu/$\sim$helton/chemjac.html} $$
for many of the   examples in this section 
(including all of those that fall under (b)),
as well as  a demonstration notebook that
readers may edit to run their own  examples. This software works
well when the number of species is small; for larger numbers,
the determinant expansion has too many terms to be handled readily.

We conclude this  preamble with some intuition
underlying the case when there are anomalous signs.
In general, we expect that for very small values of the
parameters appearing in the reaction rate functions
for a conservative reaction network
(and, in the limit, for vanishing parameter values),
the dynamics of the system augmented by inflows and outflows
will be dominated by the inflow and outflow terms,
and $\det({\partial r}/{\partial c})$ will not vanish;
moreover, if the  inflow and (linear) outflow terms dominate the
dynamics, then the equilibrium will be unique.
Examples \ref{ex:CF1} and  \ref{ex:mmCFv}
illustrate how this observation can be made rigorous
and can be used together with
Theorem \ref{thm:mm} and the
proof of Theorem \ref{thm:dissip} to conclude the existence and uniqueness of
an equilibrium for a subset of the parameter space,
even if the result does not hold for the entire
parameter space.

\subsection{The determinant expansion, its signs and uses}

\begin{exam}
\label{ex:CF1}
Consider the mass-action kinetics system given by the
chemical reaction network (\ref{rn1}),
which is an irreversible version of the network
shown in Table 1.1(i) of \cite{CF05}
(see Table \ref{table1}(i) below):
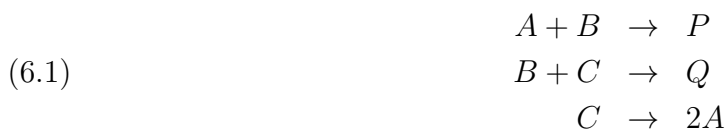
\begin{eqnarray} \label{rn1}
A+B & \to & P  \nonumber \\
B+C & \to & Q   \\ \nonumber
C   & \to & 2A
\end{eqnarray}
If we add inflow and outflow reactions for all species,
the associated dynamical system model for (\ref{rn1})
is
\begin{eqnarray} \label{syst-rn1}
&& \dot{c}_A = k_{0\to A} - k_{A\to 0}c_A - k_{A+B\to P}c_Ac_B +
2k_{C\to 2A}c_C \nonumber \\
&& \dot{c}_B = k_{0\to B} - k_{B\to 0}c_B - k_{A+B\to P}c_Ac_B -
k_{B+C\to Q}c_Bc_C  \nonumber \\
&& \dot{c}_C = k_{0\to C} - k_{C\to 0}c_C - k_{B+C\to Q}c_Bc_C
-k_{C\to 2A}c_C\\ \nonumber
&& \dot{c}_P = k_{0\to P} - k_{P\to 0}c_P + k_{A+B\to P}c_Ac_B   \\ \nonumber
&& \dot{c}_Q = k_{0\to Q} - k_{Q\to 0}c_Q + k_{B+C\to Q}c_Bc_C.
\end{eqnarray}
According to Remark 4.3 in \cite{CF05} the dynamical system above
does have multiple positive equilibria for some values of the
reaction rate parameters.

If we assume that all outflow rate constants $k_{A\to 0},...,k_{Q\to 0}$ are
equal to 1, then the determinant of the Jacobian of the reaction rate
function is:
\begin{eqnarray}
      \det({\partial r}/{\partial c}) \
      &=&
 \ -1 \ - \  k_{A+B \to P}c_A \ - \  k_{B+C\to Q}c_C
- \ k_{B+C\to Q}c_B  \nonumber \\
&& \quad
- \ k_{B+C\to Q}k_{A+B\to P}c_Ac_B  \ - \ k_{C\to 2A}
\ - \ k_{C\to 2A} k_{A+B\to P}c_A   \nonumber   \\
&&\quad -\ k_{C\to 2A}k_{B+C\to Q}c_C  \ - \  k_{A+B\to P}c_B  \ - \
k_{A+B\to P}k_{C\to 2A}c_B    \label{det-rn1}     \\
\nonumber
&& \quad -\ k_{A+B\to P}k_{B+C\to Q}c_B^2  \ - \   k_{A+B\to
P}k_{B+C\to Q}c_Bc_C \\
   \nonumber
&&\quad  + \ k_{A+B\to P}k_{B+C\to Q}k_{C\to 2A}c_Bc_C.
\end{eqnarray}
Note that there is only one positive monomial in the expansion
in (\ref{det-rn1}). The concentrations in it are $c_B c_C$,
but there is also a negative monomial with concentrations $c_B c_C$,
and the two combine to give
$$
[- \   k_{A+B\to P}k_{B+C\to Q}
+ \ k_{A+B\to P}k_{B+C\to Q}k_{C\to 2A} ] c_Bc_C.
$$
Thus if $k_{C\to 2A} \le 1$,
then the positive monomial will be dominated by
a negative monomial.
Therefore, if $k_{C\to 2A} \le 1$,
then $\det({\partial r}/{\partial c}) \neq 0$ for this network,
everywhere on ${\mathbb R}_{>0}^5$.

Note that $(m_A, m_B, m_C, m_P, m_Q)=(1,1,2,2,3)$
is a conserved
mass vector for the reaction network (\ref{rn1}).
It follows from Theorem \ref{thm:mm}
that (\ref{syst-rn1}), the dynamical system for
(\ref{rn1}), augmented with inflows and outflows (with outflow rate constants
equal to one),
has a unique positive equilibrium for all positive values of the
reaction rates such that $k_{C\to 2A} \le 1$.
Note that this uniqueness conclusion
would not follow directly from the 
theory of \cite{CF05, CF06} nor from that in \cite{BDB07}, since these
works pertain only when the determinant has the same sign
for \it all \rm rate constants and species concentrations.

The same method can be applied to conclude
   that the reversible version of the reaction network (\ref{rn1}),
augmented with inflows and outflows (with  outflow rate constants set 
equal to one),
   also has a unique positive equilibrium for all positive values of
   the reaction rates
   such that $k_{C\to 2A} \le 1$; moreover, even if the (positive)
   outflow rate constants $k_{A\to 0},...,k_{Q\to 0}$ are
   not necessarily equal to 1, the same conclusion holds
   if we know that $k_{C\to 2A} \le k_{C\to 0}$.
\qed
\end{exam}

\begin{exam}
\label{ex:CFall}
Here we summarize several examples with mass-action kinetics (in the
next subsection we consider some of these reactions with more
general kinetics).
Of the eight examples in \cite{CF05, CF06},
which are chemical reaction networks augmented with inflows and outflows
(with outflow rate constants equal to one),
two have the property that the coefficients
of the terms in
their Jacobian determinant expansion all have the same sign,
 and the other six  have all but one sign the same.
The first observation is from \cite{CF05} and the second observation, emphasizing
that there is \it only one \rm anomalous sign, is new here.
An analysis as in Example \ref{ex:CF1} can be  applied  here.
Table \ref{table1} is a list of the examples showing how many
``anomalous"
signs each determinant expansion has.
\begin{table}[htb]

\footnotesize

\begin{center}

\begin{tabular}{ccc|c}

\hline

       \rule{0pt}{10pt}
\ \ & \qquad Reaction  & &  Num. of ``anomalous" signed  \\

\ \ & \qquad network   & & terms  in  det expansion           \\

\hline

         \rule{0pt}{10pt}
(i)
& $A+B \rightleftharpoons P$
& $B+C
\rightleftharpoons Q$
& 1 \\
& $C   \rightleftharpoons 2A$   &  \\

\hline

        \rule{0pt}{10pt}
(ii)
& $A+B \rightleftharpoons P $
& $B+C
\rightleftharpoons Q $
& 0 \\

& $C+D \rightleftharpoons R $
& $D   \rightleftharpoons 2A$  & \\

\hline

        \rule{0pt}{10pt}
(iii)
& $A+B \rightleftharpoons P $
& $B+C
\rightleftharpoons Q $    \\
& $C+D \rightleftharpoons R $
& $D+E \rightleftharpoons S $
& 1
\\
& $E \rightleftharpoons 2A$  &  \\

\hline

        \rule{0pt}{10pt}
(iv)
& $A+B \rightleftharpoons P $
& $B+C
\rightleftharpoons Q $ &  0 \\
& $C   \rightleftharpoons A $ &    \\

\hline

       \rule{0pt}{10pt}
(v)
& $A+B \rightleftharpoons F $
& $A+C \rightleftharpoons G $ &   1 \\
& $C+D \rightleftharpoons B $
& $C+E \rightleftharpoons D $ &   \\

\hline

      \rule{0pt}{10pt}(vi)  & $A+B \rightleftharpoons 2A $ & & 1  \\

\hline

      \rule{0pt}{10pt}(vii) & $2A+B\rightleftharpoons 3A $ &  & 1 \\

\hline

        \rule{0pt}{10pt}(viii)& $A+2B\rightleftharpoons 3A $ & &  1  \\

\hline





\end{tabular}

\end{center}

\caption{\label{table1}Some examples of reaction networks
and the signs of coefficients  in their Jacobian determinant expansion
when augmented with inflows and
outflows (with outflow rate constants equal to one).}

\end{table}

A similar  accounting holds for examples of reaction networks
in \cite{CTF06}, see Table 1, page 8699.
These reactions involve enzymes which \cite{CTF06}
treat with mass-action kinetic models.
They find
 reaction networks 1,2,3,5,7,9 in this table
do not have any anomalous signs.
Here we point out that the remaining reaction networks,
4, 6 and 8 have very few anomalous signs.
The reaction network 4 is
$$
S + E \rightleftharpoons ES \to E + P, \quad
I + E \rightleftharpoons EI, \quad
I + ES \rightleftharpoons ESI \rightleftharpoons EI + S
$$
and has only 1 ``anomalous" sign, and the reaction network 6 is
$$
S1 + E \rightleftharpoons ES1, \quad
S2 + E \rightleftharpoons ES2, \quad
S2 + ES1 \rightleftharpoons ES1S2 \rightleftharpoons S1 + ES2,
\quad
ES1S2 \to E + P
$$
and has only 2 ``anomalous" signs.
Reaction network 8 has 4 anomalous signs  out of a total of  over
3000
terms.
Here  all reactions are augmented by  outflows 
with outflow rate constants set to one. 
(For the cases of no anomalous signs these outflow rate constants 
can be taken to be arbitrarily
small without changing the answer\footnote{See \cite{CF06iee}
for how one can eliminate outflows for the enzymes.}.)
The theory of Sections \ref{sec:mass} and \ref{sec:dynam}
applies,
if there are (arbitrarily small) inflows and outflows,
to yield that there is a  unique positive equilibrium,
for reaction networks 1,2,3,5,7,9.
It  also leaves open the possibility that one can apply the technique
in Example \ref{ex:CF1} to get a unique
positive equilibrium for certain rate constants
in reaction networks 4, 6, 8.
These applications of our theory require  finding a  conserved ``mass"
for the system without inflow and outflow,
which is easy to do in all cases.
\end{exam}

\subsection{General reaction rate functions}
In this subsection,
we follow the setup in \S \ref{sec:dynam} and
move beyond mass-action kinetics
to a very general classes of rate functions.

\begin{definition}
\label{def:incrCI}
We say that
a reaction rate function
$K_{y \to y'}$ is {\it consumptively increasing},
if for each species $i$ belonging to the support of
the vector $y$, the partial
derivative of the reaction rate function,
$\partial K_{y\to y'}/  \partial c_i$,
is
strictly positive
on  the open orthant.
\end{definition}

It is very common to assume that the reaction rate functions
$K_{y \to y'}$ are {\it consumptively increasing},
since this simply means that the rate of a reaction increases
whenever the concentration of a consumed species is increased
unilaterally.
In particular, the consumptively increasing property is true
for many common
chemical reaction rate laws, such as
many Michaelis-Menten and Hill
laws, as well as for all mass-action kinetics  \cite{keener-sneyd}.
All examples in this section entail
consumptively increasing reaction rates.

The consumptively increasing property can fail to hold for
some classes of reactions including
those involving inhibitory enzymes
and for those in which a Michaelis-Menten rate depends on
the products of the reaction \cite{Rec81}.
However, the next more lenient assumption handles
these and many additional biochemical situations.
\begin{definition}
\label{def:incrSM}
We say that
a reaction rate function
$K_{y \to y'}$ is
{\it  strictly monotone},
if for each species $i$ on which the function
$K_{y \to y'}$ actually depends,
the partial
derivative of the reaction rate function,
$\partial K_{y\to y'}/  \partial c_i$,
has
one strict sign on  the open orthant.
\end{definition}

More generally,
the main technique used in this section is to
compute the determinant expansion of 
the Jacobian
as a sum of terms,
each of which is  a product of partial derivatives
$\partial K_{y\to y'}/ \partial c_i$ where $i$ belongs
to the support of $y$.
For strictly monotone rate functions we can assign
a $\pm$  to each term according to whether
 that term is everywhere positive or negative
on the domain  $\Rnp$.
That is, strict monotonicity
 guarantees  the technique of tracking anomalous
 signs in the determinant expansion applies.


Examples \ref{ex:AS1}, \ref{mapk} and \ref{exAS3}
which involve inhibitory feedback can be written in the form
(\ref{diff_chem}) with strictly monotone rate functions.
As we observed in \S \ref{sec:exSont}
the determinant of the Jacobian, $(\frac{\pt r}{\pt c})$,
is positive for all of these situations.
However, at this point the terminology is
in place  so that we can mention the more
refined property that each of these examples has
no anomalous signs.

\begin{exam}
\label{ex:mmCFii}
   Consider the chemical reaction network (\ref{rn2}), which is the
   reversible
network shown in (ii) in Table  \ref{table1}
but, unlike in \cite{CF05}, in this example
we don't assume that the kinetics is mass-action.
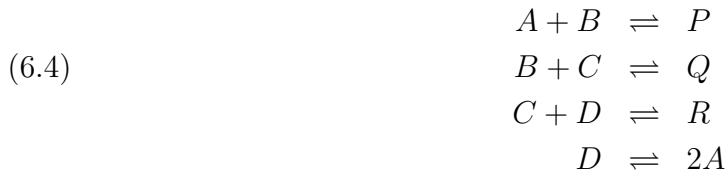
\begin{eqnarray} \label{rn2}
A+B & \rightleftharpoons & P    \nonumber \\
B+C & \rightleftharpoons & Q   \\ \nonumber
C+D & \rightleftharpoons & R   \\ \nonumber
D   & \rightleftharpoons & 2A
\end{eqnarray}
We augment this reaction with inflows and outflows
where the outflow  matrix $\lout $ is normalized to be the identity.
We suppose that each of the reaction rate functions
$K_{y \to y'}$ is
{\it consumptively increasing} as in Definition \ref{def:incrCI}.
We compute
the expansion of
$\det({\partial r}/{\partial c})$
in terms of  the partial derivatives
$\partial K_{y\to y'}(c)/ \partial c_i$,
for $i $ belonging to the support of $y$.
It is a sum of coefficients times monomials in these partial
derivatives;  the
 set of coefficients is  shown
in (\ref{det-coeffs-rn2}).
\begin{eqnarray} \label{det-coeffs-rn2}
&& \{ - 1, - 1, - 1, - 1, - 1, - 1, - 1, - 1, - 1, - 1, - 1, - 1, -1,
-1,   \nonumber             \\ \nonumber
&& -1, -1, -1, -1, -1, -1, -1, -2, -2, -2, -2, -2, -2,  -2,
\\ \nonumber
&& -2, -2, -2, -2, -2, -2,  -3, -1, -1, -1, -1, -1, -1,-1,
\\ \nonumber
&& -3, -1, -1, -1, -1, -1, -1, -1, -1, -1, -1, -1, -1, -1,
\\
&& -2, -2, -2, -2, -2, -1, -1, -1, -1, -1, -1, -1, -1,  -1,
\\ \nonumber
&& -1, -1, -1, -1, -1, -1, -1, -1, -2, -2, -2, -2, -2, -2,
\\ \nonumber
&& -1, -1, -1, -1, -1, -1, -1, -1, -2, -2, -1, -1, -1,  -1,
\\ \nonumber
&& -1, -1, -1, -1, -1, -1, -1, -2, -2, -2, -2, -2, -2, -2,
\\ \nonumber
&& -2, -1, -1, -1, -1, -1, -1, -1, -1, -2, -2, -2, -1, -1,
\\ \nonumber
&& -1, -1, -1, -1, -2, -2, -2, -1, -1, -1, -2, -1 \}
\end{eqnarray}
%

To summarize this list there are 138 terms in the expansion of
$\det({\partial r}/{\partial c})$  and the set of coefficients of
these terms contains exactly:
96 minus ones, 40 minus  twos, 2 minus  threes and no
positive terms.
This is more information than we need, since the point is
that all these numbers are negative.
This implies that $\det({\partial r}/{\partial c}) \neq 0$ for this
network, for all values of $c \in \RR^7_{>0}$ and
all ``consumptively increasing" rate laws.
Also, note that $(m_A, m_B, m_C, m_D, m_P,
m_Q, m_R) = (1,1,1,2,2,2,3) $
is a conserved mass vector for the
reaction network (\ref{rn2}).
Therefore,
we can apply Theorem \ref{thm:mm}
to
conclude that for the reaction network (\ref{rn2}) with
rate laws which are  all consumptively increasing,
and with  inflows and  outflows (with rate constants equal to one), has a unique positive
equilibrium.
\qed
\end{exam}

\begin{exam}
\label{ex:mmCFv}
This example is exactly parallel to Example \ref{ex:mmCFii}
except that now we consider
the chemical reaction network which is the
network shown in (v) in Table  \ref{table1}.
Unlike in \cite{CF05}, in this example
we don't assume that the kinetics is mass-action.
We  assume that the reaction rate functions are consumptively
increasing.
Also  we augment with inflows and outflows where
the outflow matrix   $\lout $
is normalized to be the
identity matrix.
We find that   there are 167 terms in
the expansion of
$\det({\partial r}/{\partial c})$ involving
  the partial derivatives
of $\pt K_{y \to y'}/ \pt c_i$ for
$i$ belonging to the support of $y$
 and the set of coefficients of
these terms
contains exactly: \ 146 minus ones, 20 minus  twos and one
positive term.
The positive term is
$$
K_{B \to C +D}'(c_B)\;
K_{D \to C +E}'(c_D)\; K_{A + C  \to G}^{(0, 1)}(c_A, c_C)\;
       K_{A + B  \to F}^{(1, 0)}(c_A, c_B)
$$
Here $F^{(1,0)} $ (resp. $F^{(0,1)}) $
   denotes the partial derivative of $F$
with respect to its first (respectively second) variable.

One
can think of many conditions on the reaction rate functions that
make
  the determinant have one sign on the open orthant.
Typically, the more complicated they look, the more lenient
is the assumption.
Here are some  examples.
We first note that
the reaction network (v) of Table \ref{table1}
is mass-conserving with mass vector
$m= (1, 3, 1, 2, 1, 4, 2)$.

\ben
\item
We can collect all terms containing
$K_{D \to C +E}'(c_D)\; K_{A + C  \to G}^{(0, 1)}(c_A, c_C)\;
       K_{A + B  \to F}^{(1, 0)}(c_A, c_B) $, this yields
  $(- 1+ K_{B \to C +D}'(c_B)  ) \;K_{D \to C +E}'(c_D)\;
  K_{A + C  \to G}^{(0, 1)}(c_A, c_C)\;
       K_{A + B  \to F}^{(1, 0)}(c_A, c_B) $.
Thus if we assume that $1 \geq K_{B \to C +D}'(c_B)$
for all $c_B>0$, then  the determinant is negative on the
open orthant.

\smallskip
\item
Alternatively, we can
collect all terms containing
$ K_{A + C  \to G}^{(0, 1)}(c_A, c_C)\;
       K_{A + B  \to F}^{(1, 0)}(c_A, c_B) $
and extract its coefficient which yields
\beq
\label{eq:kisgood}
- [1- K_{B \to C +D}'(c_B) ] \;K_{D \to C +E}'(c_D)
- [1 + K_{C+D \to B}^{(0, 1)} (c_C, c_D) ] \;
[ 1 + K_{C+D \to B}^{(0, 1)} (c_C, c_D)].
\eeq
Thus assuming this is negative for all positive $c_B, c_C, c_D$
makes the determinant negative on the open orthant.
We see that the second requirement is less stringent than the first.

\smallskip
\item
For any  $M> m \cdot c_{in}$ the boundary of the  set
${\Omega}_M:=\{ c\in \bRplus : m\cdot c < M\} $
contains no equilibria,
as in the proof of Theorem \ref{thm:dissip}.
A yet weaker   assumption than (1) and (2) is
that
the function (\ref{eq:kisgood}) is negative  on $\Omega_M$
for a particular $M> m\cdot c_{in}$.
\een
In cases (1) and (2),  we can apply Theorem \ref{thm:mm}
to conclude that for the reaction network in
Table \ref{table1} (v) with
rate laws which are  all consumptively increasing,
after augmentation with  inflows and  outflows (where $\lout $ is the
identity matrix), there exists a unique positive
equilibrium.
For case (3) we can apply the proof of Theorem \ref{thm:dissip}
to conclude under the same conditions that there is
a unique equilibrium in  $\Omega_M$ .
\qed
\end{exam}

\begin{rem}
We emphasize that the
   homotopy-based methods described in this paper not only imply uniqueness,
   but also \emph{existence} of a positive equilibrium for many
dynamical systems
   derived from chemical reaction networks, while the methods
    in \cite{CF05, CF06} only imply uniqueness of an equilibrium.
    Also, as we saw, our methods may be used for models containing
    very general chemical kinetics laws (not necessarily mass-action kinetics).
\end{rem}

\section{Acknowledgements}
\label{sec:acknow}

The authors wish to thank Karl Fredrickson for doing
many of the computer calculations used here and
for discussions.
Also we thank Raul Gomez for help with computations.

G. Craciun thanks the NSF
and the DOE BACTER Institute for their support.
J. W. Helton thanks the NSF and the Ford Motor company
for their support.
R. J. Williams thanks the NSF for support under
grant DMS 06-04537.


\end{document}